\documentclass[12pt,a4paper]{article}
\usepackage{amsmath}
\usepackage{amssymb}
\usepackage{amsthm}
\usepackage{amsfonts}
\usepackage{latexsym}

\theoremstyle{plain}
\newtheorem{thm}{Theorem}
\newtheorem{cor}{Corollary}
\newtheorem{lem}{Lemma}
\newtheorem{prop}{Proposition}
\theoremstyle{definition}
\newtheorem{defn}{Definition}

\theoremstyle{remark}

\newtheorem{rem}{Remark}[section]

\title{Equivariant volumes for linearized actions}
\author{Alberto Della Vedova and Roberto Paoletti
\footnote{\noindent{\bf Address.} Dipartimento di Matematica e
Applicazioni, Universit\`a degli Studi di Milano Bicocca, Via R.
Cozzi 53, 20125 Milano, Italy; {\bf e-mail}:
alberto.dellavedova@unimib.it, roberto.paoletti@unimib.it }}
\date{}
\begin{document}
\maketitle
\section{Introduction}
Let $M$ be an $n$-dimensional complex projective manifold, $\tilde
G$ a $g$-dimensional reductive connected complex Lie group, and
$\nu :\tilde G\times M\rightarrow M$ a holomorphic action. A
$\tilde G$-line bundle $(L,\tilde \nu)$ on $M$ will mean the
assignment of a holomorphic line bundle $L$ on $M$ together with a
lifting (a linearization) $\tilde \nu :\tilde G\times L\rightarrow
L$ of the action of $\tilde G$ (to simplify notation, we shall
generally leave $\tilde \nu$ understood, and denote a $\tilde
G$-line bundle by $A,\,B,\,L,\ldots$).

If $L$ is a $\tilde G$-line bundle, there is for every integer
$k\ge 0$ an induced representation of $\tilde G$ on the complex
vector space of holomorphic sections $H^0(M,L^{\otimes k})$. This
implies a $\tilde G$-equivariant direct sum decomposition
$$H^0(M,L^{\otimes k})\, \cong \,\bigoplus _{\mu \in \Lambda _+}
H^0_\mu(M,L^{\otimes k}),$$ where $\Lambda _+$ is the set of
dominant weights for a given choice of a maximal torus $\tilde
T\subseteq \tilde G$ and of a fundamental Weyl chamber. For every
dominant weight $\mu$, let us denote the associated irreducible
finite dimensional representation of $\tilde G$ by $V_\mu$. For
every $\mu\in \Lambda _+$, the summand $H^0_\mu(M,L^{\otimes k})$
is $\tilde G$-equivariantly isomorphic to a direct sum of copies
of $V_\mu$.

Suppose that the line bundle $L$ is ample. We shall address in
this paper the asymptotic growth of the dimensions
$$h^0_\mu(M,L^{\otimes k})=:\dim \left ( H^0_\mu(M,L^{\otimes
k})\right )$$ of the equivariant spaces of sections
$H^0_\mu(M,L^{\otimes k})$, for \textit{$\mu$ fixed } and
$k\rightarrow +\infty$. This problem has been the object of much
attention over the years, and has been approached both
algebraically \cite{b}, \cite{bd} and symplectically - in the
latter sense it is part of the broad and general picture revolving
around the \textit{quantization commutes with reduction} principle
\cite{gs-gq}, \cite{ggk}, \cite{ms}, \cite{s}.

In this paper, we shall study this problem under the general
assumption that the stable locus of $L$ is non-empty, $M^s(L)\neq
\emptyset$. Thus the most general result in our setting is now the
Riemann-Roch type formula for multiplicities proved in \cite{ms}.
This is a deep and fundamental Theorem - with a rather complex
symplectic proof. However, we adopt a different, more
algebro-geometric approach, taking as point of departure the
Guillemin-Sternberg conjecture for regular actions (in the sense
of \S 2 below). Our motivation was partly to understand the
leading asymptotics for singular actions by fairly elementary
algebro-geometric arguments. Besides the hypothesis that
$M^s(L)\neq \emptyset$, our arguments need an additional technical
assumption, namely that the pair $(M,L)$ admits a Kirwan
resolution with certain mildness properties. Roughly, every
\textit{divisorial} component of the unstable locus upstairs
should map to the unstable locus downstairs (Definition
\ref{defn:mild}); such resolutions will be called \textit{mild}.

One can produce examples showing that $h^0_\mu(M,L^{\otimes k})$
may not be described, in general, by an asymptotic expansion, even
if the GIT quotient $M/\!/\tilde G$ is nonsingular (but see \S 2
below and the discussion in \cite{pao-mm}). Inspired by the notion
of volume of a big line bundle \cite{del}, we shall then introduce
and study the $\mu$-equivariant volume of a $\tilde G$-line bundle
$L$, defined as
\begin{equation}
\upsilon _\mu (L)\,=:\, \limsup _{k\rightarrow +\infty}
\frac{(n-g)!}{k^{n-g}}\, h^0_\mu (X,L^{\otimes
k}).\label{eqn:equiv-volume}\end{equation}

Because it leads to a concise and simple statement, we shall focus
on the special case where the stabilizer $K\subseteq \tilde G$ of
a general $p\in M$ is a (necessarily finite) central subgroup. Our
methods can however be applied with no conceptual difficulty to
the case of an arbitrary principal type (the conjugacy class of
the generic stabilizer). This will involve singling out for each
$\mu$ the kernel $K_\mu\subseteq K$ of the action of $K$ on the
coadjoint orbit of $\mu$, and considering the contribution coming
from each conjugacy class of $K_\mu$.

Let us then assume that $K$ is a central subgroup. By restricting
the linearization to $K$, we obtain an induced character $\chi
_{K,L}:K\rightarrow \mathbb{C}^*$.

Another character $\tilde \mu_K:K\rightarrow \mathbb{C}^*$ is
associated to the choice of a dominant weight $\mu \in \lambda
_+\subseteq \frak{t}^*$. Namely, $K$ lies in the chosen maximal
torus, and $\tilde \mu _K$ is the restriction to $K$ of the
character $\tilde \mu :\tilde T\rightarrow \mathbb{C}^*$ induced
by $\mu$ by exponentiation, $\exp _{\tilde G}(\xi)\mapsto e^{2\pi
i<\mu,\xi>}$. We may define $\tilde \mu _K=:\left .\tilde
\mu\right |_K :K\rightarrow \mathbb{C}^*$ (the restriction of
$\tilde \mu$ to $K$). An alternative description of $\tilde \mu_K$
is as follows: Let $G\subseteq \tilde G$ be a maximal compact
subgroup, $T\subseteq G$ a maximal torus, and suppose $\mu \in
\frak{t}^*$, where $\frak{t}^*$ is the Lie algebra of $T$. If
$\frak{g}$ is the Lie algebra of $G$, let
$\mathcal{O}_\mu\subseteq \frak{g}^*$ be the coadjoint orbit of
$\mu$. Since $\mu$ is an integral weight, the natural K\"{a}hler
structure on $\mathcal{O}_\mu$ is in fact a Hodge form, that is,
it represents an integral cohomology class. The associated ample
holomorphic line bundle $A_\mu \rightarrow \mathcal{O}_\mu$ is a
$G$-line bundle in a natural manner. Since the action of $K$ is
trivial on $\mathcal{O}_\mu$, the linearization induces the
character $\tilde \mu _K$ on $K$.

We then have:

\begin{thm}
\label{thm:main-general-case} Let $M$ be a complex projective
manifold, $\tilde G$ a reductive complex Lie group, $\nu :\tilde
G\times M\rightarrow M$ a holomorphic action. Suppose for
simplicity that the stabilizer of a general $p\in M$ is a central
subgroup $K\subseteq \tilde G$. Let $L$ be an ample $\tilde
G$-line bundle on $M$ such that $M^s(L)\neq \emptyset$ and
admitting a mild Kirwan resolution. Let $M_0=:M/\!/\tilde G$ be
the GIT quotient with respect to the linearization $L$. Let
$\mu\in \Lambda _+$ be a dominant weight. Let $\chi _{K,L},\tilde
\mu_K:K\rightarrow \mathbb{C}^*$ be the characters introduced
above. Then:
\begin{description}
  \item[i):] If for every $r=1,\ldots,|K|$
  we have $\chi _K^r \cdot \overline{\tilde \mu }_K\neq 1$
  (the constant character equal to $1$),
  then $H^0_\mu(M,L^{\otimes k})=0$ for every $k=1,2,\ldots$;
  \item[ii):] Assume that for some $r\in \{1,\ldots,|K|\}$ we have
  $\chi _K^r \cdot \overline{\tilde \mu }_K=1$. Then
  $$\upsilon _\mu(L)=
  \dim (V_\mu)^2\cdot \mathrm{vol}(\widehat{M}_0,\widehat{\Omega }_0)>0.$$
  Here $\widetilde{M}_0$ is an orbifold,
  $\varphi: \widetilde{M}_0\rightarrow M_0$ is a partial resolution of
  singularities, $L$ induces on $\widetilde{M}_0$
  a natural nef and big line-orbibundle $\tilde L_0$,
  with first Chern class $c_1(\tilde L_0)=\left [\widetilde{\Omega}_0\right ]$, and
  $$\mathrm{vol}(\widetilde{M}_0,\widetilde{\Omega }_0)=:
  \int _{\widetilde{M}_0}\widetilde{\Omega }_0^{\wedge (n-g)}.$$
  \end{description}

\end{thm}

\section{The case $M^s(L)=M^{ss}(L)\neq \emptyset$.}

Let us begin by considering the case where the stable and
semistable loci for the $\tilde G$-line bundle are equal and
nonempty: $M^s(L)=M^{ss}(L)\neq \emptyset$. First, as we shall
make use of the Riemann-Roch formulae for multiplicities for
regular actions conjectured by Guillemin and Sternberg, and first
proved by Meinrenken \cite{mein1}, \cite{mein2}, it is in order to
recall how these algebro-geometric hypothesis translate in
symplectic terms. Let us choose a maximal compact subgroup $G$ of
$\tilde G$. Thus $G$ is a $g$-dimensional real Lie group, and
$\tilde G$ is the complexification of $G$. Let $\frak{g}$ denote
the Lie algebra of $G$. We may without loss choose a $G$-invariant
Hermitian metric $h_L$ on $L$, such that the unique covariant
derivative on $L$ compatible with $h_L$ and the holomorphic
structure has curvature $-2\pi i\Omega$, where $\Omega$ is a
$G$-invariant Hodge form on $M$. The given structure of $G$-line
bundle of $L$, furthermore, determines (and, up to topological
obstructions, is equivalent to) a moment map $\Phi =\Phi
_L:M\rightarrow \frak{g}^*$ for the action of $G$ on the
symplectic manifold $(M,\Omega)$ \cite{gs-gq}. The hypothesis that
$M^s(L)=M^{ss}(L)\neq \emptyset$ may be restated symplectically as
follows: $0\in \frak{g}^*$ is a regular value of $\Phi$, and $\Phi
^{-1}(0)\neq \emptyset$ \cite{kir1}. In this case, $P=:\Phi
^{-1}(0)$ is a connected $G$-invariant codimension $g$ submanifold
of $M$.

Let $G$ and $\Phi$ be as above. The action of $G$ on $\Phi
^{-1}(0)$ is locally free, and the GIT quotient $M/\!/\tilde
G=M^s(L)/\tilde G$ may be identified in a natural manner with the
symplectic reduction $M_0=:\Phi ^{-1}(0)/G$, and is therefore a
K\"{a}hler orbifold. The quantizing line bundle $L$ descends to a
line orbibundle $L_0$ on $M_0$.

Similar considerations apply to symplectic reductions at coadjoint
orbits sufficiently close to the origin. If the $\tilde G$-line
bundle $L$ is replaced by its tensor power $L^{\otimes k}$, the
Hodge form and the moment map get replaced by their multiples
$k\Omega$ and $\Phi _k=:k\Phi$. Given any $\mu \in \frak{g}^*$,
there exists $k_0$ such that $\mu$ is a regular value of $\Phi _k$
if $k\ge k_0$. The relevant asymptotic information about the
multiplicity of $V_\mu$ in $H^0(M,L^{\otimes k})$ may then
determined by computing appropriate Riemann-Roch numbers on these
orbifolds \cite{kawa}, \cite{mein2}.

Let $\mathcal{O}_\mu\subseteq \mathfrak{g}^*$ be the coadjoint
orbit through $\mu$; since $\mu$ is integral, the Kirillov
symplectic form $\sigma _\mu$ is a Hodge form on the complex
projective manifold $\mathcal{O}_\mu$. By the Konstant version of
the Borel-Bott theorem, there is an ample line bundle $A_\mu$ on
$\mathcal{O}_\mu$ such that $H^0(\mathcal{O}_\mu, A_\mu)$ is the
irreducible representation of $G$ with highest weight $\mu$.

Let then $M_\mu ^{(k)}$ be the Weinstein symplectic reduction of
$M$ at $\mu$ with respect to the moment map $\Phi _{k}=k\Phi _L$
($k\gg 0$). Using the normal form description of the symplectic
and Hamiltonian structure of $(M,\Omega)$ in the neighbourhood of
the coisotropic submanifold $P=\Phi ^{-1}(0)$ \cite{mein2},
\cite{gotay}, \cite{gs}, one can verify that $M_\mu^{(k)}$ is, up
to diffeomorphism, the quotient of $P\times \mathcal{O}_\mu$ by
the product action of $G$. In other words, $M_\mu^{(k)}$ is the
fibre orbibundle on $M_0=P/G$ associated to the principal
$G$-orbibundle $q:P\rightarrow M_0$ and the $G$-space
$\mathcal{O}_\mu$ (endowed with the opposite K\"{a}hler
structure); in particular its diffeotype is independent of $k$ for
$k\gg 0$. Let $p_\mu:M_\mu^{(k)}\rightarrow M_0$ be the
projection.

Let $\theta$ be a connection $1$-form for $q$ (\cite{ggk},
Appendix B). By the shifting trick, the symplectic structure
$\Omega _\mu^{(k)}$ of the orbifold $M_\mu^{(k)}$ is obtained by
descending the closed 2-form $ k\iota ^*(\Omega)+<\mu,F(\theta)>-
\sigma _\mu$ on $P\times \mathcal{O}_\mu$ down to the quotient
(the symbols of projections are omitted for notational
simplicity). The minimal coupling term $<\mu,F(\theta)>- \sigma
_\mu$ is the curvature of the line orbibundle $R_\mu=(P\times
\overline A_\mu)/G$ on $M_\mu^{(k)}$. Thus, $ \Omega _\mu^{(k)}$
is the curvature form of the line orbibundle $p_\mu^*(L_0^{\otimes
k})\otimes R_\mu$.

Let
\begin{eqnarray}\label{eqn:tildeP}
\tilde P_\mu&=:&\{(p,\mu',g)\in P\times \mathcal{O}_\mu\times G\,:
\,
g\cdot (p,\mu')=(p,\mu')\},\nonumber \\
\tilde P_{\mu,K}&=:&P\times \mathcal{O}_\mu\times K.
\end{eqnarray}
There is a natural inclusion $\tilde P_{\mu,K}\subseteq P_\mu$.
Now let $\Sigma_\mu=:\tilde P_\mu/G$, $\Sigma _{\mu,K}=:\tilde
P_{\mu,K}/G=M_\mu^{(k)}\times K$. There is a natural orbifold
complex immersion $\Sigma_\mu\rightarrow M_\mu^{(k)}$, with
complex normal orbi-bundle $N_{\Sigma_\mu}$, and $\Sigma
_{\mu,K}\subseteq \Sigma_\mu$ is the union of the $|K|$ connected
components mapping dominantly (and isomorphically) onto
$M_\mu^{(k)}$. The orbifold multiplicity of $\Sigma _{\mu,K}$ is
constant and equal to $|K|$. Let $L_0$ be the line orbi-bundle on
$M_0$ determined by descending $L$, and let $\tilde L_0$ be its
pull-back to $\Sigma _0$. Let $r$ be the complex dimension of
$\mathcal{O}_\mu$, so that $\dim M_\mu^{(k)}=n-g+r$. After
\cite{mein1} and \cite{mein2}, the multiplicity $N^{(k)}(\mu)$ of
the irreducible representation $V_\mu$ in $H^0(M,L^{\otimes k})$
is then given by:
\begin{eqnarray}\label{eqn:asymptotic-mult-abelian}
N^{(k)}(\mu)&=&\int _{\Sigma _0
}\frac{1}{d_{\Sigma_\mu}}\frac{\mathrm{Td}(\Sigma _\mu)
\mathrm{Ch}^{\Sigma _\mu}(p_\mu^*(L_0^{\otimes k})\otimes R_\mu)}{
\mathrm{D}^{\Sigma _\mu}(N_{\Sigma_\mu} )         }\nonumber\\
&=&k^{n-g}\,\frac{1}{|K|}\,\sum _{h\in K}\chi
_{K,L}(h)^k\,\overline{\tilde
\mu_K(h)} \int _{M_{\mu}^{(k)}}\,\frac{\big (k\,c_1(L_0) +c_1(R_\mu)\big )}{(n-g+r)!}^{n-g+r}\nonumber \\
&&+O(k^{n-g-1}).
\end{eqnarray}
Now suppose that $\chi _{K,L}^k\cdot \overline{\tilde \mu_K}\not
\equiv 1$. Then the action of $K$ on $\jmath^*(L^{\otimes
k})\boxtimes \overline A_\mu$ is not trivial, where $\jmath
:P\hookrightarrow M$ is the inclusion. Therefore, the fiber of
$p_\mu^*(L_0^{\otimes k})\otimes R_\mu$ on the smooth locus of
$M_0$ is a nontrivial quotient of $\mathbb{C}$, and
$N^{(k)}(\mu)=0$ in this case. If there exists $k$ such that $\chi
_{K,L}^k\cdot \overline{\tilde \mu_K}\equiv 1$, on the other hand,
the same condition holds with $k$ replaced by $k+\ell e$, where
$e$ is the period of $\chi _{K,L}$ and $\ell \in \mathbb{Z}$ is
arbitrary. Thus $k$ may be assumed arbitrarily large. Passing to
the original K\"{a}hler structure of $\mathcal{O}_\mu$ in the
computation, and recalling that $\dim (V_\mu)=(r!)^{-1}\int
_{\mathcal{O}_\mu}\,\sigma _\mu ^r$, we easily obtain:
\begin{eqnarray}\label{eqn:asymptotic-mult-abelian}
N^{(k)}(\mu)&=&\dim (V_\mu)\, \frac{k^{n-g}}{(n-g)!}\,\int
_{M_0}\,c_1(L_0)^{\wedge (n-g)}+O(k^{n-g-1}).
\end{eqnarray}

\section{The asymptotics of equivariant volumes.}

Let $f:\tilde M\rightarrow M$ be a Kirwan desingularization of the
action \cite{kir2} . This means that $f$ is a $G$-equivariant
birational morphism, obtained as a sequence of blow-ups along
$G$-invariant smooth centers, such that for all $a\gg 0$ the ample
$G$-line bundle $B=:f^*(L^{\otimes a})(-E)$ satisfies
$M^s(B)=M^{ss}(B)\supseteq f^{-1}\left (M^s(L)\right )$. Here
$E\subseteq \tilde M$ is an effective exceptional divisor for $f$.
Clearly $\upsilon _\mu (F)=\upsilon _\mu (f^*(F))$ for every
$G$-line bundle $F$ on $M$.

\begin{defn}
\label{defn:mild} Let $M_u(B)_{\mathrm{div}}\subseteq M_u(B)$ be
the divisorial part of the unstable locus of $B$; in other words,
$M_u(B)_{\mathrm{div}}$ is the union of the irreducible components
of $M_u(B)$ having codimension one in $\tilde M$. We shall say
that the Kirwan resolution $f$ is \textit{mild} if $f\left
(M_u(B)_{\mathrm{div}}\right )\subseteq M_u(L)$.\end{defn}

We have:

\begin{thm} \label{thm:equivariant-del-new}
Let $L$ be an ample $G$-line bundle on $M$ such that $M^s(L)\neq
\emptyset$. Suppose that $f:\tilde M\rightarrow M$ is a mild
Kirwan resolution of $(M,L)$. Let $H$ be a $G$-line bundle on
$\tilde M$ such that $\chi _{K,H}=1$. Let $\mu \in \Lambda _+$ be
a dominant weight. Then for any $\epsilon
>0$ there exist arbitrarily large positive integers $m$ (how large
depending on $\epsilon$ and $\mu$) such that
$$\upsilon _\mu \left (f^*(L)^{\otimes m}\otimes H^{-1}\right )
\ge m^n\, \left (\upsilon _\mu \left (f^*(L)\right
)-\epsilon\right ).$$ More precisely, this will hold whenever
$m=1+p\,e$, where $e$ is the period of $\chi _{K,L}$ and $p\in
\mathbb{N}$, $p\gg 0$.
\end{thm}

As a corollary, we obtain the following equivariant version of
Lemma 3.5 of \cite{del} (for the case of finite group actions, see
Lemma 3 of \cite{pao-ages}).

\begin{cor} Under the same hypothesis, let $H$ be a $G$-line bundle on $M$.
Then for any $\epsilon >0$ there exist arbitrarily large integers
$m>0$ (how large depending on $\epsilon$ and $\mu$) such that
$\upsilon _\mu (L^{\otimes m}\otimes H^{-1})\ge m^n\, (\upsilon
_\mu (L)-\epsilon)$. \label{cor:equivariant-del}
\end{cor}

\noindent By definition, $\upsilon _\mu(L)\ge \upsilon _\mu \left
(L^{\otimes m}\right )/m^{n-g}$ for every $m>0$. Thus Corollary
\ref{cor:equivariant-del} with $H=\mathcal{O}_M$ implies:

\begin{cor} Let $L$ be a $G$-line bundle with $M^s(L)\neq \emptyset$. Then
$$\upsilon _\mu (L)=
\limsup _{m\rightarrow +\infty}\frac{\upsilon _\mu \left
(L^{\otimes m}\right )}{m^{n-g}} .$$
\label{cor:lim-sup-equivariant-vol}
\end{cor}

Similarly,

\begin{cor} Under the same hypothesis,
$$\upsilon _\mu (L)=
\limsup _{m\rightarrow +\infty}\frac{\upsilon _\mu \left
(f^*(L^{\otimes m})(-E)\right )}{m^{n-g}} .$$
\label{cor:exceptional}
\end{cor}

\noindent \textit{Proof of Theorem \ref{thm:equivariant-del-new}.}
The proof is inspired by arguments in \cite{del}. If $\upsilon
_\mu(L)=0$, there is nothing to prove; we shall assume from now on
that $\upsilon _\mu(L)>0$, and for simplicity write $L$ for
$f^*(L)$. Thus there exists $0\le r<e$ such that $\chi
_{K,L}^r\cdot \tilde \mu _K=1$, where $e$ is the period of $\chi
_{K,L}$. If $\ell\gg 0$, by the above $H^0(\tilde M,B^{\otimes
e\ell})^G\neq 0$, so that $\upsilon _\mu \left (L^{\otimes
m}\otimes H^{-1}\right )\ge \upsilon _\mu \left (L^{\otimes
m}\otimes H^{-1}\otimes B^{-\otimes \ell e}\right )$. Thus there
is no loss of generality in replacing $H$ by $H\otimes B^{\otimes
\ell e}$ for some $\ell \gg 0$. In view of the hypothesis on $H$,
we may thus assume without loss the existence of $0\neq \sigma \in
H^0(\tilde M,H)^G\neq \{0\}$, with invariant divisor $D\in \left
|H^0(\tilde M,H)^G\right |$.

Since furthermore the class of $B$ in the $G$-ample cone
introduced in \cite {dh} lies in the interior of some chamber, the
class of $H\otimes B^{\otimes \ell e}$ lies in the same chamber
for $\ell \gg 0$. Hence we may as well assume that $H$ is a very
ample $G$-line bundle satisfying $M^s(H)=M^s(B)=M^{ss}(H)\neq
\emptyset$.

By definition of $\upsilon _\mu$, there exists a sequence $s_\nu
\uparrow +\infty$ such that
\begin{equation}\label{eqn:definition-of-v-mu}
h^0_\mu (M,L^{\otimes s_\nu})\ge \frac{s_\nu ^{n-g}}{(n-g)!}\left
(\upsilon _\mu (L)-\frac \epsilon 3\right ).\end{equation}
Necessarily $s_\nu \equiv r$ (mod $e$), $\forall \, \nu\gg 0$. We
shall show that the stated inequality holds if $p\gg 0$ and
$m=:1+p\,e$.

\begin{lem} \label{lem:subsequnce}
Fix $p\gg 0$ and let $m=:1+p\,e$. There is a sequence $k_\nu
\uparrow +\infty$ such that
\begin{equation}\label{eqn:definition-of-v-mu-adjusted} h^0_\mu
(M,L^{\otimes k_\nu})\ge \frac{k_\nu ^{n-g}}{(n-g)!}\left
(\upsilon _\mu (L)-\frac \epsilon 2\right ),\end{equation} with
$k_\nu \equiv r$ (mod $e$) and furthermore $k_\nu \equiv e$ (mod
$m$) for every $\nu$.
\end{lem}

\noindent \textit{Proof.} Let $x>0$ be an integer. We may assume
that there is a non-zero section $0\neq \tau \in H^0(M,L^{\otimes
xpe})^G$. Thus, there are injections
$$H^0_\mu (M,L^{\otimes s_\nu})\hookrightarrow
H^0_\mu (M,L^{\otimes (s_\nu+xpe)}),$$ and for $\nu \gg 0$ we have
\begin{equation}\label{eqn:definition-of-v-mu-adjusted}
h^0_\mu (X,L^{\otimes (s_\nu+xpe)})\ge \frac{s_\nu
^{n-g}}{(n-g)!}\left (\upsilon _\mu (L)-\frac \epsilon 3\right
)\ge \frac{(s_\nu +xpe)^{n-g}}{(n-g)!}\left (\upsilon _\mu
(L)-\frac \epsilon 2\right ).\end{equation} Perhaps after passing
to a subsequence, we may assume without loss of generality that
$s_\nu \equiv r'$ (mod. $m$), for a fixed $0\le r'\le m-1$. Thus,
$s_\nu=\ell _\nu m+r'$. Let now $x>0$ be an integer of the form
$x=sm+r'-e$, $s\in \mathbb{N}$. Then $xpe=x(m-1)$ and
$$s_\nu+xpe=(\ell _\nu +x)m+r'-x=(\ell _\nu +x-s)m+e.$$
Now we need only set $k_\nu=:s_\nu+xpe$.

\bigskip

Set $\ell _\nu =\left [ \frac{k_\nu}{m}\right ]$. Thus $k_\nu=\ell
_\nu m+e$.

\begin{lem}\label{lem:non-vanishing}
There exists $a>0$ such that $H^0(\tilde M,H^{\otimes m a
e}\otimes L^{\otimes -e})^G\neq \{0\}$ for every $m\ge 1$.
\end{lem}

\noindent \textit{Proof.} If $r\gg 0$, the equivalence class of
the $G$-line bundles $H^{\otimes re}\otimes L^{-e}$ in the
$G$-ample cone lie in the interior of the same chamber as the
class of $H$. Thus, they share the same stable and semistable
loci, and determine the same GIT quotient $\tilde M_0$. The
$G$-line bundles $H'=:H^{\otimes e}$ and $L'=:L^{\otimes -e}$
descend to genuine line bundles $H'_0$ and $L'_0$ on $M_0$, and
$H'$ is ample. Thus, for some $a\gg 0$ the ample line bundles
$H^{\prime \otimes a e}$ and $H^{\prime \otimes a e}\otimes L'$
are globally generated. Arguing as in \cite{gs-gq},
$H^0(M_0,H^{\prime \otimes m a}\otimes L')$ lift to $G$-invariant
sections of $H^{\otimes m a e}\otimes L^{\otimes -e}$.

\bigskip

If $R$ and $S$ are $G$-line bundles on $\tilde M$, any $0\neq
\sigma \in H^0(\tilde M,R)^G$ induces injections $H^0_\mu(\tilde
M,S)\stackrel{\otimes \sigma}{\longrightarrow} H^0_\mu(\tilde
M,R\otimes S)$. Thus, if $H^0(\tilde M,R)^G\neq 0$ then
$h^0_\mu(\tilde M,S)\le h^0_\mu(\tilde M,R\otimes S)$ for every
$\mu$. Now, in additive notation, $\ell _\nu(mL-R)=k_\nu L-(e
L+\ell_\nu R)$ for any $G$-line bundle $R$; therefore, by Lemma
\ref{lem:non-vanishing}
\begin{equation}\label{eqn:basic-inequality-kmu}
h^0_\mu \left (\tilde M,\mathcal{O}_{\tilde M}(\ell _\nu
(mL-H)\right )\ge h^0_\mu \left (\tilde M,\mathcal{O}_{\tilde
M}(k_\nu L-(\ell _\nu+am)H\right ).
\end{equation}
We may now use the $G$-equivariant exact sequences of sheaves:
$$
0\rightarrow \mathcal{O}_{\tilde M}(k_\nu L-(j+1)H)\rightarrow
\mathcal{O}_{\tilde M}(k_\nu L-jH)\rightarrow
\mathcal{O}_{D}(k_\nu L-jH)\rightarrow 0,$$ for $0\le j<s$, to
conclude inductively that
\begin{eqnarray}\label{eqn:recursive-inequality}
      h^0_\mu(\tilde M,L^{\otimes k_\nu}\otimes H^{\otimes (-s)} ) & \ge &
    h^0_\mu  (\tilde M,L^{\otimes k_\nu} )\nonumber \\
&&- \sum _{0\le j<s}h^0_\mu (D,\left .L^{\otimes k_\nu}\otimes
H^{\otimes (-s)}\right |_D )
\end{eqnarray}
We may decompose $D$ as $D=D_u+D_s$, where $D_u,D_s\ge 0$ are
effective divisors on $\tilde M$, $D_u$ is supported on the
unstable locus of $B$, $\tilde M^u(B)\subseteq M$, and no
irreducible component of $D_s$ is supported on $M^u(B)$. Let
$D_u=\sum _jD_{uj}$ and $D_s=\sum _iD_{si}$ be the decomposition
in irreducilbe components.
\begin{lem}
If $D\in \left |H^0(\tilde M,H)^G\right |$ is general, then $D_s$
is reduced, and it is nonsingular away from the unstable locus
$\tilde M^u(H)=M^u(B)$.\label{lem:D-reduced}
\end{lem}

\noindent \textit{Proof.} Perhaps after replacing $H$ by some
appropriate power we may assume that the linear series $\left
|H^0(\tilde M,H)^G\right |$ is base point free away from $\tilde
M^u(H)$. The claim then follows from Bertini's Theorem.

\bigskip
We now make use of the mildness assumption on $f$. If $m\gg 0$ is
fixed, $\nu \gg 0$ and $0\le s\le \ell _\nu+am$, then the moment
map of the (not necessarily ample) line bundle $L^{\otimes
k_\nu}\otimes H^{\otimes -s}$ is bounded away from $0$ in the
neighbourhood of $D_u$. An adaptation of the arguments in \S 5 of
\cite{gs-gq} (applied on some resolution of singularities of each
$D_{uj}$) then shows that $h^0_\mu (D_u,\left .L^{\otimes
k_\nu}\otimes H^{\otimes (-s)}\right |_{D_u} )=0$ (if $D_u$ is not
reduced, we need only filter $\mathcal{O}_{D_{uj}}(k_\nu\,L-s\,
H)$ by an appropriate chain of line bundles).

Since furthermore on each $D_{sj}$ we may find a non-vanishing
invariant section of $H$, we obtain:
\begin{eqnarray}  h^0_\mu(\tilde M,L^{\otimes k_\nu}\otimes H^{\otimes (-s)} ) & \ge &
    h^0_\mu  (\tilde M,L^{\otimes k_\nu} ) \\
&&- \sum _{0\le j<s}h^0_\mu (D_s,L^{\otimes k_\nu}\otimes
H^{\otimes (-s)}\otimes
\mathcal{O}_{D_{s}} ) \nonumber \\
 & \ge &h^0_\mu
(\tilde M,L^{\otimes k_\nu} )- sh^0_\mu (D_s,L^{\otimes
k_\nu}\otimes
\mathcal{O}_{D_{s}} )\nonumber \\
& \ge &h^0_\mu (\tilde M,L^{\otimes k_\nu} )- sh^0_\mu
(D_s,(L\otimes B)^{\otimes k_\nu}\otimes \mathcal{O}_{D_{s}}
)\nonumber
\end{eqnarray}

\begin{prop}
There exists $C>0$ constant such that if $D\in \left |H^0(\tilde
M,H)^G\right |$ is general then
$$
h^0_\mu(D_s,(L\otimes B)^{\otimes k}\otimes
\mathcal{O}_{D_{s}})\le Ck^{n-g-1}$$ for every $k\gg 0$.
\label{prop:divisor-bound}
\end{prop}

\noindent \textit{Proof.} Given the equivariant injective morphism
of structure sheaves $\mathcal{O}_{D_s}\longrightarrow \bigoplus
_i\mathcal{O}_{D_{si}},$ we may as well assume that $D_s$ is a
reduced and irreducible $G$-invariant divisor, descending to a
Cartier divisor $D_0$ on the quotient.

Let $R=:L\otimes B$, with associated moment map $\Phi _R$. By the
generality in its choice, we may assume that $D_s\in \left
|H^0(\tilde M,H)^G\right |$ is non-singular in the neighbourhood
of $\Phi _R^{-1}(0)$, and is transversal to it. In fact, the
singular locus of $\left |H^0(\tilde M,H)^G\right |$ is the
unstable locus of $H$. Furthermore, by the arguments of Lemma 3 in
\cite{pao-mm} and compactness one can see the following: There
exist a finite number of holomorphic embeddings  $\varphi
_i:B\rightarrow \tilde M$, where $B\subseteq \mathbb{C}^{n-g}$ is
the unit ball, satisfying i): $\varphi _i(B)\subseteq \Phi
_R^{-1}(0)$; ii): as a submanifold of $\Phi _R^{-1}(0)$, $\varphi
_i(B)$ is transversal to every $G$-orbit; iii): the union $\bigcup
_i \varphi _i(B)$ maps surjectively onto $\tilde M/\!/\tilde G$.
In view the local analytic proof of Bertini's theorem in
\cite{gh}, we may assume that $D$ is transversal to each $\varphi
_i(B)$. By $G$-invariance, it is then transversal to all of $\Phi
_R^{-1}(0)$.

Let $g:\tilde D_s\rightarrow D_s$ be a $G$-equivariant resolution
of singularities \cite{eh}, \cite{ev}. For $s\gg 0$,
$g^*(R^{\otimes s})(-F)$ is an ample $G$-line bundle on $\tilde
D_s$, where $F$ is some effective exceptional divisor. Since $0\in
\frak{g}^*$ is a regular value of $g\circ \Phi _R:\tilde
D_s\rightarrow \frak{g}^*$ and belongs to its image, the same
holds for the moment map of $g^*(R^{\otimes s})(-F)$, for $s\gg
0$. Having fixed $s\gg 0$, let us also choose $r\gg 0$ such that
$H^0(\tilde M,g^*(R^{\otimes sre})(-reF))^G\neq 0$. The choice of
$0\neq \sigma \in H^0(\tilde D_s,g^*(R^{\otimes sre})(-reF))^G$
determines injections
$$
H^0_\mu(D_s,R^{\otimes k}\otimes
\mathcal{O}_{D_{s}})\,\stackrel{\otimes \sigma ^{\otimes
k}}{\longrightarrow} \, H^0_\mu(\tilde D_s,R^{\otimes
k(1+sre)}(-rekF)\otimes \mathcal{O}_{\tilde D_{s}}).$$ This
implies the statement by the arguments in \S 2, since $\dim
(\tilde D_s)=n-1$.

\bigskip

Given (\ref{eqn:definition-of-v-mu}),
(\ref{eqn:recursive-inequality}) and Proposition
\ref{prop:divisor-bound}, we get
\begin{equation}\label{eqn:final-recursive-ineq}
h^0_\mu(\tilde M,L^{\otimes k_\nu}\otimes H^{\otimes (-s)} ) \ge
\frac{k_\nu^{n-g}}{(n-g)!}\left (\upsilon
_\mu(L)-\frac{\epsilon}{2}\right )- s\,C\, k_\nu ^{n-g-1}.
\end{equation}
Now, in view of (\ref{eqn:basic-inequality-kmu}), we set $s=\ell
_\mu +am$ to obtain:
\begin{equation}
\begin{array}{lll}

    h^0_\mu \left ( \tilde M, \mathcal{O}_{\tilde M} (
    \ell _\mu (m L-H) )\right ) & \ge &
    \frac{ k_\nu ^{n-g} }{ (n-g)! }
    \left ( \upsilon _\mu (L)-\frac{\epsilon}{2} \right )-C
    (\ell _\nu+am) k_\nu ^{n-g-1} \\
      & \ge &\frac{\ell _\nu ^{n-g} m^{n-g}}{(n-g)!}\left (\upsilon _\mu(L)-\frac{\epsilon}{2}\right )\\
       & & -
     C(\ell _\nu+am)(\ell _\nu+1)^{n-g-1}m^{n-g-1}.
  \end{array}
\end{equation}
The proof of Theorem \ref{thm:equivariant-del-new} follows by
taking $\ell _\nu\gg m\gg 1$ (see \cite{del}, Lemma 3.5).

\begin{rem}
The arguments used in the proof of Theorem
\ref{thm:equivariant-del-new} may be applicable in other
situations. For example, suppose that $L$ is a $G$-ample line
bundle with $\mathrm{vol}_0(L)>0$; assume that the equivalence
class of $L$ in the $G$-ample cone \cite{dh} lies on a face of
measure zero, and that - say - the $G$-ample line bundles in the
interior of an adjacent chamber have unstable locus of codimension
$\ge 2$. If $A$ is a $G$-ample line bundle in the interior of the
chamber, one may apply the previous arguments to tensor powers of
the form $L^{\otimes k}\otimes A$.\end{rem}

\section{Proof of Theorem \ref{thm:main-general-case}.}

Given the Kirwan resolution $f:\tilde M\rightarrow M$, for $m\gg
0$ the equivalence classes of the ample $G$-line bundles $f^*\left
(L^{\otimes m}\right )(-E)$ in the $G$-ample cone of $\tilde M$
all lie to the interior of the same chamber. Therefore, they
determine the same GIT quotient $\tilde M_0=\tilde M/\!/\tilde G$.
The latter is an $(n-g)$-dimensional complex projective orbifold,
which partially resolves the singularities of $M/\!/\tilde G$
\cite{kir2}. Being $G$-line bundles on $\tilde M$, $f^*(L)$ and
$\mathcal{O}_{\tilde M}(-E)$ descend to line orbi-bundles on
$\tilde M_0$. Fixing $G$-invariant forms $\tilde \Omega$ and
$\Omega _{-E}$ on $\tilde M$ representing the first Chern class of
$f^*(L)$ and $\mathcal{O}_{\tilde M}(-E)$, we obtain forms $\tilde
\Omega _0$ and $\Omega _{-E\,0}$ on $\tilde M_0$.

By Corollary \ref{cor:exceptional}, we have
$$\upsilon _\mu(L)=\limsup _{m\rightarrow +\infty}
\frac{\upsilon _\mu\left(f^*\left (L^{\otimes m}\right )(-E)\right
)}{m^{n-g}}.$$ By the results in \S 2, under the appropriate
numerical hypothesis,
$$\upsilon _\mu\left(f^*\left (L^{\otimes m}\right )(-E)\right )=
\dim (V_\mu)^2\cdot \mathrm{vol} \left (\tilde M_0,m\tilde \Omega
_0+\Omega _{-E\,0}\right ).$$ The statement follows.

\end{document}